\RequirePackage[leqno]{amsmath}
\documentclass[leqno]{amsart}
\usepackage{amsmath,mathtools}
\makeatletter
\usepackage[utf8]{inputenc}     % Accept different input encodings
\usepackage{multicol,url}
\usepackage[inline]{enumitem} 
\usepackage{color, graphics}
\usepackage[demo]{graphicx}
\usepackage{efbox,graphicx}
\efboxsetup{linecolor=gray,linewidth=0.8pt}
\usepackage{floatrow}
\usepackage{float}
\usepackage{pgfplots}
\usepackage{latexsym, amsmath, amssymb}
\usepackage{caption}
\usepackage{subcaption}
\usepackage{pgfplots}
\usepackage{amsfonts}
\usepackage{pifont}
\usepackage{enumitem}
\setlist[itemize]{topsep=0pt,after=\vspace{1.5\baselineskip}}
\usepackage[colorlinks=true]{hyperref}
\usepackage{empheq}
\usepackage[T1]{fontenc}
\usepackage{tabularx}
\usepackage[leftcaption]{sidecap}
\sidecaptionvpos{figure}{m}
\usepackage{tikz}
\usepackage{caption}
\usepackage{tikz-cd}
\usetikzlibrary{decorations.markings}
\usetikzlibrary{decorations.pathmorphing}
\usepackage{tabularx,multirow,array} 
\usepackage[margin=0.5in]{geometry}
\usetikzlibrary{backgrounds,automata}

\pgfplotsset{compat=1.18}
\usepackage{xparse}
\ExplSyntaxOn
\int_new:N \g_tohi_const_int
\int_new:N \g_tohi_const_sub_int
\tl_new:N  \g_tohi_const_char_tl

\cs_new_protected:Nn \tohi_print_constant:nn 
{
  #1 \textsubscript {#2} 
} 

\NewDocumentCommand\resetconstants{m}
{
 \int_gincr:N \g_tohi_const_int
 \int_gzero:N \g_tohi_const_sub_int
 \tl_gset:Nn  \g_tohi_const_char_tl {#1}
}

\NewDocumentCommand\const{m}
{
  \tl_if_exist:cTF
   {
    c_tohi_const_\int_use:N\g_tohi_const_int _#1_tl
   }
   {
    \tl_use:c {c_tohi_const_\int_use:N\g_tohi_const_int _#1_tl }
   }
   {
    \int_gincr:N \g_tohi_const_sub_int
    \tl_const:cx {c_tohi_const_\int_use:N\g_tohi_const_int _#1_tl }
     { \exp_not:N\tohi_print_constant:nn {\g_tohi_const_char_tl }{\int_use:N \g_tohi_const_sub_int}}
    \tl_use:c {c_tohi_const_\int_use:N\g_tohi_const_int _#1_tl }
   }
}

\ExplSyntaxOff
\newcommand{\inlineitem}[1][]{%
\ifnum\enit@type=\tw@
    {\descriptionlabel{#1}}
  \hspace{\labelsep}%
\else
  \ifnum\enit@type=\z@
       \refstepcounter{\@listctr}\fi
    \quad\@itemlabel\hspace{\labelsep}%
\fi}
\DeclarePairedDelimiter\abs{\lvert}{\rvert}
\DeclarePairedDelimiter\norm{\lVert}{\rVert}
\DeclarePairedDelimiter\tonda{(}{)}
\DeclarePairedDelimiter\quadra{[}{]}
\DeclarePairedDelimiter\graffa{\{}{\}}
\newcommand{\into}{\int_\Omega}

\setlist[itemize]{noitemsep, topsep=0pt}

\def\R{\mathbb R} \def\N{\mathbb N} 
\def\tmax{\ensuremath{T_{\textup{max}}}}

\def\R{\mathbb R} \def\N{\mathbb N} 

\def\nvq{\abs{\nabla v}^2}
 
\def\dqvq{\abs*{D^2 v}^2}
\def\absnv{\abs{\nabla v}}
\def\nivz{\norm{v_0}_{L^{\infty}(\Omega)}}
\newcommand{\norminf}[1]{\norm*{#1}_{L^{\infty}(\Omega)}}
\def
\@cite
#1#2{[#1\if@tempswa, #2\fi]}
\makeatother
\newtheorem{theorem}{Theorem}[section]
\newtheorem{corollary}[theorem]{Corollary}

\newtheorem{lemma}[theorem]{Lemma}

\newtheorem{remark}[theorem]{Remark}

\title[Boundedness criteria for a chemotaxis consumption model with gradient nonlinearities] 
{Boundedness criteria for a chemotaxis consumption model with gradient nonlinearities}
\author[Alessandro Columbu]{Alessandro Columbu}
\keywords{Chemotaxis, Attraction, Consumption, Boundedness. \\
\textit{$^\star$Corresponding author}: alessandro.columbu2@unica.it}

\begin{document}

\maketitle
{
\medskip
\centerline{Dipartimento di Matematica e Informatica}
\centerline{Universit\`{a} di Cagliari}
\centerline{Via Ospedale 72, 09124. Cagliari (Italy)}
\medskip
}
\bigskip
%\tableofcontents
\resetconstants{c}

\begin{abstract}
    This work deals with the consumption chemotaxis problem
    \begin{equation*}
    \begin{cases*}
        u_t = \Delta u - \chi \nabla \cdot u\nabla v + \lambda u - \mu u^2 - c \lvert \nabla u \rvert^\gamma, & \text{in $\Omega\times(0,\tmax)$}, \\
        v_t = \Delta v - uv, & \text{in $\Omega\times(0,\tmax)$}, 
    \end{cases*}
\end{equation*}
in a bounded and smooth domain $\Omega\subset\R^n$, $n\geq 3$, under Neumann boundary conditions, for $\chi,\lambda,\mu,c>0$, $\tmax\in(0,\infty]$ and for $u_0,v_0$ positive initial data with a certain regularity. We will show that the problem has a unique and uniformly bounded classical solution for $\gamma\in\bigl(\frac{2n}{n+1},2\bigr]$. Moreover, we have the same result for $\gamma=\frac{2n}{n+1}$ and a condition that involves the parameters $c,\mu,n,\chi$ and the initial data.
\end{abstract}

\section{Introduction and main result}
Chemotaxis is a phenomenon in which bacteria or, in general, uni- or multi-cellular organisms, direct their movement according to the presence of certain chemicals in the environment, that usually is impenetrable. In \cite{Keller-1971-MCabbr} we can find one of the first tests to write a model for this phenomenon, which, in the general case, is described by two coupled PDEs with homogeneous Neumann boundary conditions:
\begin{equation*}\label{general}
    \begin{cases}
        u_t = \nabla\cdot \tonda*{A(u,v)\nabla u - B(u,v)\nabla v} + h(u) & \text{in $\Omega\times(0,\tmax)$}, \\
        \tau v_t = \Delta v + E(u,v), & \text{in $\Omega\times(0,\tmax)$}, \\
        u_\nu = v_\nu = 0, & \text{on $\partial\Omega\times(0,\tmax)$}, \\
        u(x,0) = u_0(x), \; \tau v(x,0) = \tau v_0(x), & x \in \Omega,
    \end{cases}
\end{equation*}
where $\Omega\subset\R^n$ is a bounded and smooth domain and $\tmax\in(0,\infty])$ is the maximum existence time of the solution. The first equation gives us information on how the distribution of the biological material changes over time, conditioned by the natural diffusion with law $A(u,v)$, by the signal sensitivity $B(u,v)$, and a source $h(u)$. The second equation, that can be elliptic or parabolic for $\tau=0$ or $\tau=1$ respectively, describes the evolution of the chemical signal which affects the biological material: in particular, this chemical can be produced by the cells or, like the case we want to study, it can be consumed, with the law $E(u,v)$. Our specific case deals with linear diffusion $A(u,v)=1$, a sensitivity $B(u,v)=\chi u$ and the logistic source $h(u)=\lambda u-\mu u^2-c\lvert \nabla u \rvert^\gamma$. For the second equation, since we are studying a consumption model, $E(u,v)=-uv$, different from the production model that is written $E(u,v)=-v+f(u)$, where $f(u)$ is a positive function. The two main concerns with this type of problem is to find conditions that ensure the global boundedness of the solution or the so called \emph{blow-up}, an uncontrolled aggregation of the biological material. For the linear production case, the reader can see for instance \cite{WinklAggreabbr}, researchers found a condition which involves the mass $m=\into u$ and the chemosensitivity coefficient $\chi$. Moreover, the blow-up situation is very much studied in the signal-production version and the literature is full of these results. On the contrary, for the consumption model, there is nothing on that.
\subsection{The consumption model}In our case, except when $n\in\{1,2\}$ where the solutions are always bounded, results concerning global boundedness involve both the chemosensitivity $\chi$ and the maximum of the initial condition $\nivz$. The first of this results is contained in \cite{TaoBounabbr}, in which the author found $\chi\nivz<\frac{1}{6(n+1)}$ for $n\geq 3$. This has been, later, improved in \cite{BaghaeiKhelghati-nolog} by $\chi\nivz<\frac{\pi}{\sqrt{2(n+1)}}$. By adding a dampening logistic term, one can find conditions even if it is not verified what we saw previously. One example is given by \cite{LankeitWangConsumptLogisticabbr}, where 
\begin{equation*}
    u_t = \Delta u - \chi \nabla \cdot u\nabla v + \lambda u - \mu u^2, \quad \text{in $\Omega\times(0,\tmax)$},
\end{equation*}
and the solution is globally bounded if $\mu$ is larger than a quantity which increases with $\chi\nivz$. Another similar result, where $\mu$ has to be large enough, can be found in \cite{MarrasViglialoroMathNachabbr}, where the authors consider both nonlinear diffusion and sensitivity.
\subsection{The logistic in the form of \texorpdfstring{$\boldsymbol{h(u,\nabla u)}$}{h(u)}} According to \cite{Souplet_Gradientabbr}, the evolution of a certain species in an habitat is influenced, not only by births and natural deaths, also by the accidental deaths, which can be described with the term $-\lvert \nabla u \rvert^\gamma$. In the literature, we can mention here \cite{chipot_weissler,filar,kawohl_peletier}, there are some results in which gradient nonlinearities are studied in order to ensure or prevent blow-up formations for a semilinear heat equation with Dirichlet conditions. Moreover, in \cite{IshidaLankeitVigliloro-Gradientabbr} the authors studied a production chemotaxis problem with a gradient dependent logistic and they achieved, for the boundedness of the solution, a condition on the parameter $\gamma$.
\subsection{Main result}
The idea is to consider a consumption model with a dampening logistic that depends also on $\nabla u$:
\begin{equation}\label{consum}
    \begin{cases}
        u_t = \Delta u - \chi \nabla \cdot u\nabla v + \lambda u - \mu u^2 - c \lvert \nabla u \rvert^\gamma, & \text{in $\Omega\times(0,\tmax)$}, \\
        v_t = \Delta v - uv, & \text{in $\Omega\times(0,\tmax)$}, \\
        u_\nu = v_\nu = 0, & \text{on $\partial\Omega\times(0,\tmax)$}, \\
        u(x,0) = u_0(x), \; v(x,0) = v_0(x), & x \in \Omega,
    \end{cases}
\end{equation}

% \textrm{ for some } \delta\in(0,1) \textrm{ and } n\in \N, \Omega \subset \R^n \textrm{ is a bounded domain of class } C^{2+\delta},\\
%     u_0, v_0, w_0: \bar{\Omega}  \rightarrow \R^+ , \textrm{ with } u_0,  v_0,   w_0 \in C_\nu^{2+\delta}(\bar\Omega)=\{\psi \in C^{2+\delta}(\bar{\Omega}): \psi_\nu=0 \textrm{ on }\partial \Omega\}, 

where $\Omega\subset\R^n$ is a bounded domain of class $C^{2+\delta}$, with $\delta\in(0,1)$ and the initial conditions are such that
\begin{equation}\label{reginitialconditions}
    u_0,v_0\colon \bar{\Omega}  \rightarrow \R^+, \text{ with } u_0,  v_0 \in C_\nu^{2+\delta}(\bar\Omega)=\{\psi \in C^{2+\delta}(\bar{\Omega}): \psi_\nu=0 \textrm{ on }\partial \Omega\}.
\end{equation}
For $p>1$ and $\eta>0$, we define the functions
    \begin{equation}\label{kappa1}
        K_1(p,n)\coloneqq\frac{p}{p+1}\tonda*{\frac{8(4p^2+n)}{p(p+1)}}^{\frac 1p}\tonda*{\frac{p(p-1)}{2}}^{\frac{p+1}{p}},
    \end{equation}
    and
    \begin{equation}\label{kappa2}
        K_2(p,n,\eta)\coloneqq\frac{2p^{\frac{p+1}{2}}(p+n+\eta-1)^{\frac{p+1}{2}}}{p+1}\tonda*{\frac{8(4p^2+n)(p-1)}{p(p+1)}}^{\frac{p-1}{2}}.
    \end{equation}
Moreover, let $M \coloneqq \max{\graffa*{\into u_0, \frac \lambda\mu \abs{\Omega}}}$. We enunciate now the main result of this paper:
    % \begin{equation}\label{condsommacmu}
    %          c\frac n2 \tonda*{\frac{2n}{n^2+3n-2}}^{\frac{2n}{n+1}}\mathcal{C}M^{-\frac{2}{n+1}}+\mu\frac n2>K_1\tonda*{\frac n2,n}+K_2\tonda*{\frac n2,n,0}
    % \end{equation}
\begin{theorem}\label{teo}
    Let $\Omega\subset\R^n$, $n\geq 3$ be a bounded domain of class $C^{2+\delta}$, with $\delta\in(0,1)$, $u_0,v_0$ as in \eqref{reginitialconditions}, $\lambda,\mu,c,\chi>0$. There exists $\mathcal{C}>0$ such that whenever
    \begin{enumerate}[label={($\mathcal{A}_{\arabic*}$)}]
        \item \label{a1} $\displaystyle \frac{2n}{n+1}<\gamma\leq 2$;
    \end{enumerate}
    or
    \begin{enumerate}[resume, label={($\mathcal{A}_{\arabic*}$)}]
        \item \label{a2} $\displaystyle \gamma=\frac{2n}{2+1}$ \qquad and \qquad $\displaystyle c\frac n2 \tonda*{\frac{2n}{n^2+3n-2}}^{\frac{2n}{n+1}}\mathcal{C}M^{-\frac{2}{n+1}}+\mu\frac n2>K_1\tonda*{\frac n2,n}\nivz^{\frac 4n}\chi^{\frac 2n(n+2)} +K_2\tonda*{\frac n2,n,0}\nivz^n$, 
    \end{enumerate}
    problem \eqref{consum} admits a unique and uniformly bounded classical solution.
\end{theorem}
\begin{remark}
Let us dedicate to the case $\gamma=\frac{2n}{n+1}$ and rewrite now the inequality in condition \ref{a2} as
\begin{equation}\label{eqespl}
    c\mathcal{F}M^{-\frac{2}{n+1}}+ \frac n2 \mu > \mathcal{K},
\end{equation}
where 
\begin{equation*}
    \mathcal{F}=\frac n2 \tonda*{\frac{2n}{n^2+3n-2}}^{\frac{2n}{n+1}}\mathcal{C} \quad \text{and} \quad \mathcal{K}=K_1\tonda*{\frac n2,n}\nivz^{\frac 4n}\chi^{\frac 2n(n+2)} +K_2\tonda*{\frac n2,n,0}\nivz^n.
\end{equation*} 
We can observe what follows:
% [$\blacktriangleright$]
\begin{enumerate}[label={($\mathcal{B}_{\arabic*}$)}]
    \item Fix $c=0$. We can write \eqref{eqespl} as $\mu = \Bar{\mu} > \frac 2n \mathcal{K}$, and this is another way to see the result contained in \cite{LankeitWangConsumptLogisticabbr}, where the authors use $p=n$, so the condition for the boundedness contained in that work is a special case of this one.
    \item   In the case $\mu < \frac 2n \mathcal{K}$ and $c>0$, we have two different possibilities:
        \begin{enumerate}
        \item [$\blacktriangleright$] if $M=\into u_0$, we obtain
        \begin{equation*}
            c\mathcal{F}\tonda*{\into u_0}^{-\frac{2}{n+1}}+ \frac n2 \mu > \mathcal{K}, \quad \text{that leads to} \quad \into u_0<\tonda*{\frac{c\mathcal{F}}{\mathcal{K}-\frac n2 \mu}}^{\frac{n+1}{2}},
        \end{equation*}
        so we have a condition on the initial mass;
        \item [$\blacktriangleright$] if $M=\frac{\lambda}{\mu}\abs{\Omega}$ we have
        \begin{equation}\label{eqmu}
            \mathcal{E}\mu^{\frac{2}{n+1}}+\frac n2 \mu > \mathcal{K} \quad \text{where} \quad \mathcal{E}=\frac{c\mathcal{F}}{(\lambda\abs{\Omega})^{\frac{2}{n+1}}},
        \end{equation}
        and so there exists $\mu<\Bar{\mu}$ such that \eqref{eqmu} is satisfied.
    \end{enumerate}
    \item In the situation where $\mu=\frac 2n \mathcal{K}$ and $c>0$, relation \eqref{eqespl} is simply true.
\end{enumerate}
\end{remark}
\section{Preliminaries}
We begin this section with a result taken from \cite[Lemma 2.2]{LankeitWangConsumptLogisticabbr}:
\begin{lemma} Let $q\in [1,\infty)$. Then for any $w\in C^2(\Bar{\Omega})$ satisfying $w\frac{\partial w}{\partial \nu}=0$ on $\partial \Omega$, the following inequality holds
    \begin{equation}\label{grad4qquadro}
        \norm{\nabla w}^{2q+2}_{L^{2q+2}(\Omega)}\leq 2(4q^2+n)\norm{w}^2_{L^\infty(\Omega)} \norm*{\abs{\nabla w}^{q-1}D^2 w}^2_{L^2(\Omega)}.
    \end{equation}
\end{lemma}
In addition, the next Lemma provides some technical inequalities which will be used later on.
\begin{lemma}\label{lemmap}
Let $n\in\N$ and the lower bound of condition \ref{a1} be fixed. Then for every $p>1$ and
\begin{equation*}
    \theta(p,\gamma)=\frac{\frac{p+\gamma-1}{\gamma}\tonda*{1-\frac{1}{p+1}}}{\frac{p+\gamma-1}{\gamma}+\frac 1n-\frac 1\gamma} \quad \sigma(p,\gamma)=\frac{\gamma(p+1)}{p+\gamma-1} \quad \Bar{\theta}(p)=\frac{\frac p2-\frac 12}{\frac p2-\frac 12 +\frac 1n} \quad \hat{\theta}=\theta\tonda*{p,\frac{2n}{n+1}} \quad \hat{\sigma}=\sigma\tonda*{p,\frac{2n}{n+1}},
\end{equation*}
these inequalities hold
\begin{table}[H]
 \begin{subequations}
     \begin{subtable}[h]{0.24\textwidth}
        \centering
        \begin{equation}\label{theta}
        0<\theta<1,
        \end{equation}
   \end{subtable}
   \hfill
    \begin{subtable}[h]{0.24\textwidth}
        \centering
        \begin{equation}\label{thetasigma}
        0<\frac{\theta\sigma}{\gamma}<1,
        \end{equation}
    \end{subtable}
    \hfill
    \begin{subtable}[h]{0.24\textwidth}
        \centering
        \begin{equation}\label{thetabar}
        0<\Bar{\theta}<1,
        \end{equation}
    \end{subtable}
    \hfill
    \begin{subtable}[h]{0.24\textwidth}
        \centering
        \begin{equation}\label{thetasigmahat}
        \frac{\hat{\sigma}\hat{\theta}(n+1)}{2n}=1.
        \end{equation}
    \end{subtable}
\end{subequations}
\end{table}
\begin{proof}
 First of all, we can see that \eqref{theta},\eqref{thetabar} and the lower bound of \eqref{thetasigma} are true for every $p>1$. We need $\gamma>\frac{2n}{n+1}$ in order to prove the upper bound of \eqref{thetasigma}.
\end{proof}
\end{lemma}
\begin{lemma}\label{lemmal}
    Let $y$ be a positive real number verifying $y\leq k(y^l+1)$ for some $k > 0$ and $0 < l < 1$. Then $y\leq \max{\graffa*{1,(2k)^{\frac{1}{1-l}}}}$.
\begin{proof}
    For $y<1$ is trivial. If $y\geq 1$ then $y^l\geq 1$ and, subsequently, $y\leq k(y^l+1)\leq 2ky^l$. So we have the thesis.
\end{proof}
\end{lemma}
\section{Local existence and extension criterion}
The proofs of the next three results are an adaption to the consumption case of the works \cite{IshidaLankeitVigliloro-Gradientabbr,ACVgradient}.
\begin{lemma}\label{existencet}Let $u_0,v_0$ satisfy \eqref{reginitialconditions}, $\lambda,\mu,c,\chi>0$, $\gamma\geq 1$. Then there exists $T>0$ and a unique solution $(u,v)$ of system \eqref{consum} in $\Omega\times(0,T)$ such that
\begin{equation*}
    (u,v)\in \tonda*{C^{2+\delta,1+\frac{\delta}{2}}( \Bar{\Omega} \times [0, T])}^2. 
\end{equation*}
\begin{proof}
The main idea is articulated as follows (for the details see \cite{ACVgradient}). First, for any $R>0$ and for some $0<T\leq1$, we define the closed convex set 
\begin{equation*}
    S_{T}=\{0\le u \in C^{1,\frac{\delta}{2}}(\bar{\Omega}\times [0,T]):  \|u(\cdot,t)-u_0\|_{L^\infty(\Omega)}\le R,\, \; \text{for\ all}\ t\in[0,T]\},
\end{equation*}
and we choose one element $\tilde{u}$ from it. From the second equation, we find $v$ with a certain regularity that allows us to find the unique solution $u$ to the first equation. In particular, $u\in S_{T}$ and this implies that the compact operator  $\Phi(\tilde{u})=u$ is such that $\Phi(S_{T})\subset S_{T}$. So, by the Schauder's fixed point theorem, let $u$ be the fixed point of $\Phi$. Now, by taking $T$ as the initial time and $u(\cdot,T)$ as the initial condition, we can find a solution $(\hat{u},\hat{v})$ defined in $\Omega\times[T,\hat{T}]$ for some $\hat{T}>0$, that is the prolongation of $(u,v)$. We can repeat this procedure in order to obtain a maximum interval time of existence $[0,\tmax)$.
\end{proof}
\end{lemma}

\begin{corollary}\label{existencetmax}
Let $u_0,v_0$ be as in \eqref{reginitialconditions}, $\lambda,\mu,c,\chi>0$ and $\gamma\geq 1$. There exists $\tmax \in (0,\infty]$ such that $u,v>0$ in $\Omega\times[0,\tmax)$, and 
\begin{equation}\label{extension}
    \text{if }\tmax<\infty, \text{ then } \limsup_{t\nearrow\tmax}{\tonda*{\norm*{u(\cdot,t)}_{C^{2+\delta}(\bar{\Omega})}+\norm*{v(\cdot,t)}_{C^{2+\delta}(\bar{\Omega})}}}=\infty.
\end{equation}
\end{corollary}
\noindent As the reader may noticed, the fact that $\tmax$ is finite, does not ensure the explosion in the $L^\infty(\Omega)$-norm. In order to have this, we need $\gamma\leq2$.
\begin{lemma}
    Let $u_0,v_0$ be as in \eqref{reginitialconditions}, $\lambda,\mu,c,\chi>0$ and $\gamma\in[1,2]$. Then, if $u\in L^\infty((0,\tmax),L^\infty(\Omega))$, then $\tmax=\infty$.
\begin{proof}
For the details we refer to \cite{ACVgradient}. The main idea is to proceed by contradiction: if $\tmax$ is finite, the boundedness of $u$ in $L^{\infty}(\Omega)$ implies, thanks to $\gamma\in[1,2]$, the uniform-in-time boundedness of $\norm*{u(\cdot,t)}_{C^{2+\delta}(\bar{\Omega})}+\norm*{v(\cdot,t)}_{C^{2+\delta}(\bar{\Omega})}$ so, from \eqref{extension}, we have an inconsistency.
\end{proof}
\end{lemma}
\noindent From now on, $(u,v)$ will be the unique solution of problem \eqref{consum} defined on $(0,\tmax)$ and every $c_i$ is a positive constant.
\begin{lemma}
    Let $u_0,v_0$ be as in \eqref{reginitialconditions}, $\lambda,\mu,c,\chi>0$ and $\gamma\geq1$. Then 
    \begin{equation}\label{boundednessmass}
        \into u \leq M \coloneqq \max{\graffa*{\into u_0, \frac \lambda\mu \abs{\Omega}}} \quad \text{on $(0,\tmax)$},
    \end{equation}
    \begin{equation}\label{boundednessv}
        \norm{v(\cdot,t)}_{L^{\infty}(\Omega)}\leq \norm{v_0}_{L^{\infty}(\Omega)} \quad \text{for all $t\in (0,\tmax)$}.
    \end{equation}
\begin{proof}
By integrating the first equation and using the Neumann boundary conditions, we have
\begin{equation}\label{mass1}
    \frac{d}{dt}\into u\leq \lambda \into u -\mu \into u^2 \quad \text{for all $t\in(0,\tmax)$}.
\end{equation}
Moreover, an application of Hölder's inequality leads to 
\begin{equation*}
    \into u \leq \abs{\Omega}^{\frac 12}\tonda*{\into u^2}^{\frac 12} \quad \text{and, consequently, to} \quad -\into u^2\leq -\frac{1}{\abs{\Omega}}\tonda*{\into u}^2 \quad \text{for every $t\in(0,\tmax)$}.
\end{equation*}
By inserting the last one into \eqref{mass1} and calling $m=\into u$, we obtain this initial problem
\begin{equation*}
\begin{dcases*}
    m'\leq \lambda m - \frac{\mu}{\abs{\Omega}}m^2 & \text{on $(0,\tmax)$},\\
    m(0)=\into u_0,
\end{dcases*}
\end{equation*}
which ensures that $m\leq \max{\graffa*{\into u_0, \frac \lambda\mu \abs{\Omega}}}$ for all $t\in(0,\tmax)$. Moreover, \eqref{boundednessv} is a consequence of the parabolic maximum principle because 
\begin{equation*}
    \max_{\Omega\times(0,\tmax)} v(x,t)\leq \max_{\Bar{\Omega}}v(x,0).
\end{equation*}
\end{proof}
\end{lemma}
\begin{lemma}\label{lemmavq}
Suppose there exists $\tau\in(0,\min{\{1,\tmax\}})$ and $p\in[1,n]$ such that $\norm{u(\cdot,t)}_{L^p(\Omega)}\leq \const{lk1}$ for all $t\in(\tau,\tmax)$. Then for every $\delta\in(0,\tmax-\tau)$ and for all $q<\frac{np}{n-p}$ there exists $\const{lk2}$ such that
\begin{equation*}
    \norm{\nabla v(\cdot,t)}_{L^q(\Omega)}\leq \const{lk2} \quad \text{for every $t\in(\tau+\delta,\tmax)$}.
\end{equation*}
\begin{proof}
As in \cite[Section 2]{HorstWinkabbr}, we define 
\begin{equation*}
    A_q w\coloneqq -\Delta w \quad \text{where}\quad w\in D(A_q)=\graffa*{\psi\in W^{2,q}(\Omega)\text{ : }\frac{\partial \psi}{\partial \nu}=0 \text{ on }\partial\Omega}.
\end{equation*}
It is known that:
\begin{equation}\label{inclusioneD}
        D\tonda*{(A_q+1)^\beta}\xhookrightarrow{}W^{1,q}(\Omega) \quad \text{if }\beta>\frac 12,
\end{equation}
and, for every $1\leq p<q\leq\infty$ and $w\in L^p(\Omega)$, there exists a positive constant $\hat{c}$ such that
\begin{equation}\label{disA+1}
        \norm*{(A_q+1)^\beta e^{-tA_q}w}_{L^q(\Omega)}\leq \hat{c}t^{-\beta-\frac n2\tonda*{\frac 1p -\frac 1q}}e^{(1-\mu)t}\norm{w}_{L^p(\Omega)}.
\end{equation}
Let $q<\frac{np}{n-p}$ and $\beta>\frac 12$ such that $\beta+\frac n2\tonda*{\frac 1p-\frac 1q}<1$. Now we write the second equation as $v_t=-(A_q+1)v+v(1-u)$ and we use the representation formula to affirm that
\begin{equation*}
    v(\cdot,t)=e^{-(t-\tau)}v(\cdot,\tau)+\int_\tau^t e^{-(t-s)(A_q+1)}v(\cdot,s)(1-u(\cdot,s))\,ds \quad \text{for all $t\in[\tau,\tmax)$}.
\end{equation*}
By applying $(A_q+1)^\beta$ to both sides and using \eqref{disA+1} and \eqref{boundednessv} we have
\begin{equation*}
\begin{split}
    \norm*{(A_q+1)^\beta v(\cdot,t)}_{L^q(\Omega)}\leq & \const{lk4}(t-\tau)^{-\beta-\frac n2\tonda*{1-\frac 1q}}\norm{v(\cdot,\tau)}_{L^1(\Omega)}+\const{lk5}\int_\tau^t (t-s)^{-\beta-\frac n2\tonda*{\frac 1p -\frac 1q}}e^{-\mu(t-s)}\norm*{v(\cdot,s)(1-u(\cdot,s))}_{L^p(\Omega)}\, ds \\
    \leq & \const{lk4}\delta^{-\beta-\frac n2\tonda*{1-\frac 1q}}\norm{v(\cdot,0)}_{L^1(\Omega)}+ \const{lk5}\norm{v(\cdot,0)}_{L^\infty(\Omega)} \int_\tau^t (t-s)^{-\beta-\frac n2\tonda*{\frac 1p -\frac 1q}}e^{-\mu(t-s)}\tonda*{\abs{\Omega}^{\frac 1p}+\norm*{u(\cdot,s)}_{L^p(\Omega)}}\, ds \\
    \leq &\const{lk6} +\const{lk7}\int_0^\infty z^{-\beta-\frac n2\tonda*{\frac 1p -\frac 1q}}e^{-\mu z}\,dz \leq \const{lk8} \quad \text{on $(0,\tmax)$}.
\end{split}
\end{equation*}
Finally, \eqref{inclusioneD} gives the claim.
\end{proof}
\end{lemma}
\noindent The next Lemma transforms the boundedness of $u$ on some $L^p$-norm in the boundedness on the $L^{\infty}$-norm.
\begin{lemma}\label{lemmalimitato}
    Let $\tmax\in(0,\infty]$, $\frac n2<p<n$,  $\lambda,\mu,c,\chi>0$ and $\gamma\in[1,2]$. If $u\in L^\infty((0,\tmax),L^p(\Omega))$, then $u\in L^\infty((0,\infty),L^\infty(\Omega))$.
\begin{proof}
We write $u(\cdot,t)$ using the variation-of-constants formula for all $t\in(0,\tmax)$ and $t_0=\max{\{0,t-1\}}$
    \begin{equation}\label{riga1}
        u(\cdot,t)=e^{(t-t_0)\Delta}u(\cdot,t_0)-\int_{t_0}^t e^{(t-s)\Delta}\nabla\cdot \tonda*{u(\cdot,s)\nabla v(\cdot,s)}\,ds+\int_{t_0}^t e^{(t-s)\Delta}\tonda*{\lambda u(\cdot,s)-\mu u^2(\cdot,s)-c\abs{\nabla u(\cdot,s)}^\gamma}\,ds.
    \end{equation}
Focusing on the logistic term, we can say that
\begin{equation*}
    \lambda u-\mu u^2-c\abs{\nabla u}^\gamma\leq \lambda u-\mu u^2\leq \sup_{\xi>0}{(\lambda\xi-\mu\xi^2)}\coloneqq\const{f1},
\end{equation*}
and, in conjunction with the positivity of the heat semigroup, we have
\begin{equation}\label{riga2}
    \int_{t_0}^t e^{(t-s)\Delta}\tonda*{\lambda u(\cdot,s)-\mu u^2(\cdot,s)-c\abs{\nabla u(\cdot,s)}^\gamma}\,ds\leq \const{f1}(t-t_0)\leq \const{f1}.
\end{equation}
In addition, thanks to \cite[Lemma 1.3]{WinklAggreabbr} we know that there exist positive constants $C_f$ and $\Hat{C}_f$ such that, for all $t>0$ and $f\in L^1(\Omega)$
\begin{equation}\label{nhs1}
    \norminf{e^{t\Delta}f}\leq C_f\tonda*{1+t^{-\frac n2}}\norm{f}_{L^1(\Omega)}
\end{equation}
and for all $t>0$, $r>1$ and $f\in C^1(\Bar{\Omega})$ with $\frac{\partial f}{\partial \nu}=0$ on $\partial \Omega$
\begin{equation}\label{nhs2}
    \norminf{e^{t\Delta}\nabla \cdot \nabla f}\leq \hat{C}_f\tonda*{1+t^{-\frac 12 -\frac{n}{2r}}}e^{-\lambda_1 t}\norm*{\nabla f}_{L^r(\Omega)}
\end{equation}
where $\lambda_1$ is the first nonzero eigenvalue of $-\Delta$ in $\Omega$ under Neumann boundary conditions. Consequently, if $t\leq 1$, which is equivalent to $t_0=0$, we have for the parabolic maximum principle 
\begin{equation}\label{u01}
    \norminf{e^{(t-t_0)\Delta}u(\cdot,t_0)}\leq \norminf{u_0}.
\end{equation}
Instead, if $t>1$, that implies $t-t_0=1$, we can use \eqref{nhs1} with $f=u(\cdot,t_0)$ and the boundedness of the mass so to write
\begin{equation}\label{u02}
    \norminf{e^{(t-t_0)\Delta}u(\cdot,t_0)}\leq C_f \tonda*{1+t^{-\frac n2}}\norm{u(\cdot,t)}_{L^1(\Omega)}\leq \const{jk}.
\end{equation}

Moreover, for all $r\in(n,q)$, by applying \eqref{nhs2}, we know that
\begin{equation}\label{riga3}
    \norm*{\int_{t_0}^t e^{(t-s)\Delta} \nabla\cdot \tonda*{u(\cdot,s)\nabla v (\cdot,s)}\,ds}_{L^{\infty}(\Omega)}\leq \const{as1} \int_{t_0}^t \tonda*{1+(t-s)^{-\frac 12 -\frac{n}{2r}}}\norm*{u(\cdot,s)\nabla v(\cdot,s)}_{L^{r}(\Omega)} \,ds.
\end{equation}
We now define 
\begin{equation*}
    A(t')=\sup_{t\in(0,t')} \norm{u(\cdot,t)}_{L^{\infty}(\Omega)},
\end{equation*}
so that we can find
\begin{equation}\label{riga4}
\begin{split}
    \norm*{u(\cdot,s)\nabla v(\cdot,s)}_{L^r(\Omega)}\leq & \norm{u(\cdot,s)}_{L^{\frac{rq}{q-r}}(\Omega)} \norm*{\nabla v(\cdot,s)}_{L^q(\Omega)}\\
    \leq & \norm{u(\cdot,t)}_{L^{\infty}(\Omega)}^{1-\frac{p(q-r)}{rq}} \norm{u(\cdot,t)}_{L^{p}(\Omega)}^{\frac{p(q-r)}{rq}} \norm*{\nabla v(\cdot,s)}_{L^q(\Omega)}\\
    \leq &\const{dff} A(t')^{1-\frac{p(q-r)}{rq}}=\const{dff} A^l(t'),
\end{split}
\end{equation}
with $l<1$, where in the last step we applied Lemma \ref{lemmavq}. By recalling \eqref{riga1}, \eqref{riga2},  \eqref{u01}, \eqref{u02}, \eqref{riga3} and \eqref{riga4} we have $A(t')\leq \const{df1}\tonda*{A^l(t')+1}$, that implies, thanks to Lemma \ref{lemmal}, $\norm{u(\cdot,t)}_{L^\infty(\Omega)}\leq A(t')\leq \max{\graffa*{1,(2\const{df1})^{\frac{1}{1-l}}}}=\const{df2}$ for all $t\in(0,t')$ and, in particular, $\const{df2}$ does not depend on the time.
\end{proof}
\end{lemma}

\noindent We present now some technical inequalities which we will use in the next section and from now on, the proof will follow \cite{LankeitWangConsumptLogisticabbr}.
\begin{lemma}
Let $u_0,v_0$ be as in \eqref{reginitialconditions}, $\lambda,\mu,c,\chi>0$ and $\gamma\geq1$. Then there exists $C>0$ such that
\begin{equation}\label{boundednessgradv2}
    \into \abs{\nabla v(\cdot,t)}^2\leq C \quad \text{for all $t\in(0,\tmax)$}.
\end{equation}
\begin{proof}The reader can find this proof in \cite[Lemma 3.4]{LankeitWangConsumptLogisticabbr}.
\end{proof}
\end{lemma}

\begin{lemma}Under the hypotheses of Lemma \ref{lemmap}, we have that for every $p>1$ and for all $\varepsilon,\hat{c}>0$ there exists $\const{b2},\const{b7},\const{b5},\const{b4}$ such that these inequalities hold
    \begin{equation}\label{gnup}
        (\lambda p+1)\into u^p \leq \frac{p-1}{p}\into \abs*{\nabla u^{\frac p2}}^2 +\const{b2} \quad \text{on $(0,\tmax)$},
    \end{equation}
    \begin{equation}\label{gngradv2}
        \into \abs{\nabla v}^{2p}\leq \frac p4 \into \abs{\nabla v}^{2p-2}\abs*{D^2 v}^2+\const{b7} \quad \text{on $(0,\tmax)$},
    \end{equation}
    \begin{equation}\label{gammamaggiore}
        \hat{c}\into u^{p+1} \leq \varepsilon \into \abs*{\nabla u^{\frac{p+\gamma-1}{\gamma}}}^\gamma +\const{b5}  \quad \text{on $(0,\tmax)$},
    \end{equation}
    \begin{equation}\label{gammauguale}
        \hat{c}\into u^{p+1} \leq \hat{c}(2C_{GN})^{\hat{\sigma}} M^{\frac{2}{n+1}} \into \abs*{\nabla u^{\frac{p(n+1)+n-1}{2n}}}^{\frac{2n}{n+1}} +\const{b4} \quad \text{on $(0,\tmax)$},
    \end{equation}
    \begin{equation}\label{youngchi}
        \frac{p(p-1)}{2}\chi^2\into u^p \nvq \leq \frac p4 \into \absnv^{2p-2}\dqvq+K_1(p,n)\nivz^{\frac 2p}\chi^{\frac{2(p+1)}{p}}\into u^{p+1} \quad \text{on $(0,\tmax)$},
    \end{equation}
    \begin{equation}\label{youngv0}
        p(p+n-1+\eta)\nivz^2 \into u^2\absnv^{2p-2}\leq \frac p4 \into \absnv^{2p-2}\dqvq+K_2(p,n,\eta)\nivz^{2p}\into u^{p+1} \quad \text{on $(0,\tmax)$},
    \end{equation}
where $K_1(p,n)$ and $K_2(p,n,\eta)$ are defined in \eqref{kappa1} and \eqref{kappa2}.
\begin{proof}In order to prove \eqref{gnup} we make use of the Gagliardo--Nirenberg and Young inequalities in conjunction with \eqref{thetabar} and \eqref{boundednessmass}:
    \begin{equation*}
    \begin{split}
        (\lambda p+1)\into u^p = (\lambda p+1) \norm*{u^{\frac p2}}^2_{L^2(\Omega)}\leq & \const{b8}\tonda*{\norm*{\nabla u^{\frac p2}}^{\Bar{\theta}}_{L^2(\Omega)}\norm*{u^{\frac p2}}^{1-\Bar{\theta}}_{L^{\frac 2p}(\Omega)}+\norm*{u^{\frac p2}}_{L^{\frac 2p}(\Omega)}}^2 \\
        \leq & \const{b9} \tonda*{\into \abs*{\nabla u^{\frac p2}}^2}^{\Bar{\theta}}+\const{b8qq}
        \leq  \frac{p-1}{p} \into \abs*{\nabla u^{\frac p2}}^2 +\const{b2} \quad \text{on $(0,\tmax)$}.
    \end{split}
    \end{equation*}
For \eqref{gngradv2}, we operate in same way but recalling \eqref{boundednessgradv2}:
    \begin{equation*}
    \begin{split}
        \into \abs{\nabla v}^{2p}=&\into \tonda*{\abs{\nabla v}^p}^2=\norm*{\abs{\nabla v}^p}^2_{L^2(\Omega)}\leq \const{b10} \tonda*{\norm*{\nabla \abs*{\nabla v}^p}^{\Bar{\theta}}_{L^2(\Omega)}\norm*{\abs{\nabla v}^p}^{1-\Bar{\theta}}_{L^{\frac 2p}(\Omega)}+\norm*{\abs*{\nabla v}^p}_{L^{\frac 2p}(\Omega)}}^2\\
        \leq & \frac{1}{4p} \into \abs*{\nabla \abs*{\nabla v}^p}^2 +\const{b11} \leq \frac p4 \into \abs*{\nabla v}^{2p-2}\abs*{D^2 v}^2 +\const{b7} \quad \text{on $(0,\tmax)$}.
    \end{split}
    \end{equation*}
Again, with the same inequalities but recalling \eqref{theta}, \eqref{thetasigma} and, again, \eqref{boundednessmass}, we are able to prove relation \eqref{gammamaggiore}:
    \begin{equation*}
    \begin{split}
        \hat{c}\into u^{p+1}=&\hat{c}\norm*{u^\frac{p+\gamma-1}{\gamma}}^\sigma_{L^\sigma(\Omega)}\leq \hat{c}C_{GN}^\sigma \tonda*{\norm*{\nabla u^{\frac{p+\gamma-1}{\gamma}}}^\theta_{L^{\gamma}(\Omega)} \norm*{u^{\frac{p+\gamma-1}{\gamma}}}^{1-\theta}_{L^{\frac{\gamma}{p+\gamma-1}}(\Omega)}+\norm*{u^{\frac{p+\gamma-1}{\gamma}}}_{L^{\frac{\gamma}{p+\gamma-1}}(\Omega)}}^\sigma\\
        \leq & \const{b13} \tonda*{\into \abs*{\nabla u^{\frac{p+\gamma-1}{\gamma}}}^\gamma}^{\frac{\sigma\theta}{\gamma}}+\const{b14}\leq \varepsilon \into \abs*{\nabla u^{\frac{p+\gamma-1}{\gamma}}}^\gamma+\const{b5} \quad \text{on $(0,\tmax)$}.
    \end{split}
    \end{equation*}
    If $\gamma=\frac{2n}{n+1}$ then \eqref{thetasigmahat} is verified, so we cannot apply Young's inequality:
    \begin{equation*}
    \begin{split}
        \hat{c}\into u^{p+1}=&\hat{c}\norm*{u^\frac{p(n+1)+n-1}{2n}}^{\hat{\sigma}}_{L^{\hat{\sigma}}(\Omega)} \\ \leq & \hat{c}C_{GN}^{\hat{\sigma}}\tonda*{\norm*{\nabla u^{\frac{p(n+1)+n-1}{2n}}}^{\hat{\theta}}_{L^{\frac{2n}{n+1}}(\Omega)} \norm*{u^{\frac{p(n+1)+n-1}{2n}}}^{1-\hat{\theta}}_{L^{\frac{2n}{p(n+1)+n-1}}(\Omega)}+\norm*{u^{\frac{p(n+1)+n-1}{2n}}}_{L^{\frac{2n}{p(n+1)+n-1}}(\Omega)}}^{\hat{\sigma}}\\ 
        \leq & \hat{c}(2C_{GN})^{\hat{\sigma}} M^{\frac{2}{n+1}}\into \abs*{\nabla u^{\frac{p(n+1)+n-1}{2n}}}^{\frac{2n}{n+1}} +\const{b4} \quad \text{on $(0,\tmax)$}.
    \end{split}
    \end{equation*}
In order to obtain bound \eqref{youngchi}, firstly we use Young's inequality and then \eqref{grad4qquadro} so to have
\begin{equation*}
\begin{split}
    \frac{p(p-1)}{2}\chi^2\into u^p \nvq\leq & \varepsilon_1 \into \absnv^{2(p+1)}+\frac{p}{p+1}(\varepsilon_1(p+1))^{-\frac 1p}\quadra*{\frac{p(p-1)}{2}\chi^2}^{\frac{p+1}{p}}\into u^{p+1}\\
    \leq &\varepsilon_1 2(4p^2+n)\nivz^2 \into \absnv^{2p-2}\dqvq + \frac{p}{p+1}(\varepsilon_1(p+1))^{-\frac 1p}\quadra*{\frac{p(p-1)}{2}\chi^2}^{\frac{p+1}{p}}\into u^{p+1} \\
    = & \frac p4 \into \absnv^{2p-2}\dqvq+K_1(p,n)\nivz^{\frac 2p}\chi^{\frac{2(p+1)}{p}}\into u^{p+1} \quad \text{on $(0,\tmax)$},
\end{split}
\end{equation*}
where in the last step we choose $\varepsilon_1=\frac{p}{8(4p^2+n)\nivz^2}$. As to bound \eqref{youngv0} we have
\begin{equation*}
\begin{split}
     p(p+n-1+\eta)\nivz^2 \into u^2\absnv^{2p-2}\leq & \varepsilon_1 \into \absnv^{2(p+1)}+\frac{2}{p+1}(\varepsilon_1\frac{p+1}{p-1})^{-\frac{p-1}{2}}\quadra*{p(p+n-1+\eta)\nivz^2}^{\frac{p+1}{2}}\into u^{p+1}\\
    \leq &\varepsilon_1 2(4p^2+n)\nivz^2 \into \absnv^{2p-2}\dqvq \\
     &+ \frac{2}{p+1}(\varepsilon_1\frac{p+1}{p-1})^{-\frac{p-1}{2}}\quadra*{p(p+n-1+\eta)\nivz^2}^{\frac{p+1}{2}}\into u^{p+1} \\
    = & \frac p4 \into \absnv^{2p-2}\dqvq+K_2(p,n,\eta)\nivz^{2p}\into u^{p+1} \quad \text{on $(0,\tmax)$},
\end{split}
\end{equation*}
with the same value of $\varepsilon_1$ that we used in the previous proof.
\end{proof}
\end{lemma}

\section{Conditions for the existence of bounded classical solutions}
The goal is to reach a differential inequality of the type $\varphi(t)+\varphi'(t)\leq C$ for some positive $C$ and, to do that, we begin with the following two Lemmas, which will be used together in the proof of Lemma \ref{lemmanormap}.
\begin{lemma}Let $u_0,v_0$ be as in \eqref{reginitialconditions}, $\lambda,\mu,c,\chi>0$ and $\gamma>1$. Then, for any $p\in[1,\infty)$ the following holds
\begin{equation}
    \frac{d}{dt}\into u^p\leq -\frac{2(p-1)}{p}\into \abs*{\nabla u^{\frac p2}}^2+\frac{p(p-1)}{2}\chi^2\into u^p\abs{\nabla v}^2+\lambda p\into u^p-\mu p\into u^{p+1}-cp\tonda*{\frac{p+\gamma-1}{\gamma}}^{-\gamma}\into \abs*{\nabla u^{\frac{p+\gamma-1}{\gamma}}}^\gamma,
\end{equation}
for every $t\in(0,\tmax)$.
\begin{proof}We calculate the time derivative of the functional $\into u^p$ and we apply Young's inequality, so to have, on $(0,\tmax)$,
\begin{equation*}
\begin{split}
        \frac{d}{dt}\into u^p=&-p(p-1)\into u^{p-2}\abs{\nabla u}^2+\chi p(p-1)\into u^{p-1}\nabla u\cdot \nabla v +\lambda p \into u^p -\mu p \into u^{p+1} -c p \into u^{p-1} \abs{\nabla u}^\gamma \\
        \leq & -\frac{2(p-1)}{p}\into \abs*{\nabla u^{\frac p2}}^2 +\frac{\chi^2 (p-1)}{2}\into u^p \abs{\nabla v}^2+\lambda p \into u^p -\mu p \into u^{p+1} -c p \tonda*{\frac{p+\gamma-1}{\gamma}}^{-\gamma}\into \abs*{\nabla u^{\frac{p+\gamma-1}{\gamma}}}^\gamma. \qedhere
\end{split}
\end{equation*}
\end{proof}
\end{lemma}

\begin{lemma}
Let $u_0,v_0$ be as in \eqref{reginitialconditions}, $\lambda,\mu,c,\chi>0$ and $\gamma\geq1$. Then there exists $\eta>0$ such that, for any $p\in[1,\infty)$ this inequality is satisfied
\begin{equation}
    \frac{d}{dt}\into \absnv^{2p}+p\into \absnv^{2p-2} \dqvq \leq p(p+n-1+\eta)\nivz^2 \into u^2\absnv^{2p-2}+\const{c1} \quad \text{for all $t\in(0,\tmax)$}.
\end{equation}
\begin{proof}This proof is contained in \cite[Lemma 4.2]{LankeitWangConsumptLogisticabbr}.
\end{proof}
\end{lemma}

\begin{lemma}\label{lemmanormap} Let $u_0,v_0$ be as in \eqref{reginitialconditions} and $\lambda,\mu,c,\chi>0$.
If $\gamma>\frac{2n}{n+1}$ or $\gamma=\frac{2n}{n+1}$ and
\begin{equation}\label{condsommacmuconetaep}
    cp\tonda*{\frac{2n}{(n+1)p+n-1}}^{\frac{2n}{n+1}}(2C_{GN})^{-\hat{\sigma}}M^{-\frac{2}{n+1}}+\mu p\geq K_1(p,n)\nivz^{\frac 2p}\chi^{\frac{2(p+1)}{p}}+K_2(p,n,\eta)\nivz^{2p}
\end{equation}
then
    \begin{equation*}
        \into u^p+\into \absnv^{2p}\leq C \quad \text{for all $t\in(0,\tmax)$}.
    \end{equation*}

\begin{proof}
Let $\varphi(t)=\into u^p+\into \absnv^{2p}$, then, by using \eqref{gnup}, \eqref{gngradv2}, \eqref{youngchi} and \eqref{youngv0} we have
    \begin{equation*}
    \begin{split}
        \varphi'(t)+\varphi(t)\leq& -\frac{2(p-1)}{p}\into \abs*{\nabla u^{\frac p2}}^2+\frac{p(p-1)}{2}\chi^2\into u^p\nvq+(\lambda p+1)\into u^p-\mu p\into u^{p+1}\\ &-cp\tonda*{\frac{p+\gamma-1}{\gamma}}^{-\gamma}\into \abs*{\nabla u^{\frac{p+\gamma-1}{\gamma}}}^\gamma -p\into \absnv^{2p-2} \dqvq \\ &+p(p+n-1+\eta)\nivz^2 \into u^2\absnv^{2p-2}+\into \absnv^{2p} +\const{c1}\\
        \leq & -\frac{p-1}{p} \into \abs*{\nabla u^{\frac p2}}^2-\frac p4 \into \absnv^{2p-2} \dqvq -cp\tonda*{\frac{p+\gamma-1}{\gamma}}^{-\gamma}\into \abs*{\nabla u^{\frac{p+\gamma-1}{\gamma}}}^\gamma +\const{da1} \\ &+\tonda*{K_1(p,n)\nivz^{\frac 2p}\chi^{\frac{2(p+1)}{p}}+K_2(p,n,\eta)\nivz^{p-1}-\mu p}\into u^{p+1}\\
        \leq & \tonda*{K_1(p,n)\nivz^{\frac 2p}\chi^{\frac{2(p+1)}{p}}+K_2(p,n,\eta)\nivz^{p-1}-\mu p}\into u^{p+1}-cp\tonda*{\frac{p+\gamma-1}{\gamma}}^{-\gamma}\into \abs*{\nabla u^{\frac{p+\gamma-1}{\gamma}}}^\gamma +\const{da1}
    \end{split}
    \end{equation*}
If $\gamma>\frac{2n}{n+1}$ we can apply \eqref{gammamaggiore} with $\varepsilon$ small enough. If, instead, $\gamma=\frac{2n}{n+1}$ we apply \eqref{gammauguale} and use the condition \eqref{condsommacmuconetaep}.
% \begin{equation*}
%     \tonda*{K_1(p,n)+K_2(p,n,\eta)-\mu p}C_{GN}^{\hat{\sigma}}M^{\frac{2}{n+1}}-cp\tonda*{\frac{(n+1)p+n-1}{2n}}^{-\frac{2n}{n+1}}\leq0
% \end{equation*}
In both cases we arrive at
\begin{equation*}
    \varphi'(t)+\varphi(t)\leq C
\end{equation*}
that ensures the thesis.
\end{proof}
\end{lemma}

\subsection*{Proof of Theorem \ref{teo}}
If \ref{a1} is valid, Lemmata \ref{lemmanormap} and \ref{lemmalimitato} ensure the thesis. Instead, if \ref{a2} is satisfied, since $K_1(p,n),K_2(p,n,\eta),C_{GN}$ and $\hat{\sigma}$ are continuous functions of the vaiables $p$ and $\eta$, we can say that there exist $p>\frac n2$ and $\eta>0$ such that the condition \eqref{condsommacmuconetaep} is true, and again, both Lemmas give us the conclusion.

\bigskip 
	
\subsubsection*{Acknowledgements}
AC is member of the {\em Gruppo Nazionale per l’Analisi Matematica, la Probabilità e le loro Applicazioni} (GNAMPA) of the Istituto Nazionale di Alta Matematica (INdAM). AC is supported by
\begin{enumerate}
   \item [$\bullet$] the GNAMPA-INdAM Project \textit{Equazioni alle derivate parziali nella modellizazione di fenomeni reali} \\(CUP--E53C22001930001);
    \item [$\bullet$] the GNAMPA-INdAM Project \textit{Problemi non lineari di tipo stazionario ed evolutivo} (CUP--E53C23001670001)
    \item [$\bullet$] the research
project {\em Analysis of PDEs in connection with real
phenomena}, CUP F73C22001130007, funded by
\href{https://www.fondazionedisardegna.it/}{Fondazione di Sardegna}, annuity 2021.
\end{enumerate}

\bibliographystyle{plain}

\end{document}